\documentclass[12pt]{article}
\usepackage{amsmath,amssymb,epsfig}
\usepackage{mathrsfs}
\usepackage{color}

\newcommand{\openN}{{{\rm I}\kern-.16em {\rm N}}}
\newcommand{\openZ}{{{\rm Z}\kern-.28em{\rm Z}}}
\newcommand{\openR}{{{\rm I}\kern-.16em {\rm R}}}

\newcommand{\rr}{\mathbb{R}}

\newcommand{\hh}{\mathbb{H}}

\def\ri{{\rm i}\;\!}

\def\a{{\bf a}}
\def\b{{\bf b}}
\def\c{{\bf c}}

\def\i{{\bf i}}
\def\j{{\bf j}}
\def\k{{\bf k}}

\def\m{{\bf m}}
\def\n{{\bf n}}

\def\p{{\bf p}}

\def\u{{\bf u}}
\def\v{{\bf v}}
\def\w{{\bf w}}

\def\A{{\bf A}}
\def\B{{\bf B}}

\def\M{{\bf M}}
\def\N{{\bf N}}
\def\I{{\bf I}}

\def\q{{\bf q}}
\def\P{{\cal P}}
\def\Q{{\cal Q}}
\def\U{{\cal U}}
\def\V{{\cal V}}
\def\W{{\cal W}}
\def\BCH{\mathrm{BCH}}

\def\half{\textstyle\frac{1}{2}\displaystyle}
\def\quarter{\textstyle\frac{1}{4}\displaystyle}
\def\eighth{\textstyle\frac{1}{8}\displaystyle}

\newcommand{\be}{\begin{equation}}
\newcommand{\ee}{\end{equation}}
\newcommand{\ba}{\begin{eqnarray}}
\newcommand{\ea}{\end{eqnarray}}
\newcommand{\bi}{\begin{itemize}}
\newcommand{\ei}{\end{itemize}}

\newtheorem{dfn}{Definition}
\newtheorem{rmk}{Remark}
\newtheorem{lma}{Lemma}
\newtheorem{crl}{Corollary}
\newtheorem{prn}{Proposition}
\newtheorem{thm}{Theorem}

\newtheorem{exm}{Example}
\newcommand{\prf}{\noindent{{\bf Proof} :\ }}
\newcommand{\QED}{\vrule height 1.4ex width 1.0ex depth -.1ex\ \medskip}

\begin{document}

\title{Minkowski products of unit quaternion sets}

\author{
Rida~T.~Farouki \\
Department of Mechanical and Aerospace Engineering, \\
University of California, Davis, CA 95616, USA 
\medskip \\
Graziano Gentili\\
Dipartimento di Matematica e Informatica ``U. Dini,'' \\
Universit\`a di Firenze, Viale Morgagni 67/A, I--50134 Firenze, Italy 
\medskip \\
Hwan Pyo Moon \\
Department of Mathematics, \\ 
Dongguk University--Seoul, Seoul 04620, Republic of Korea
\medskip \\
Caterina Stoppato \\
Dipartimento di Matematica e Informatica ``U. Dini,'' \\
Universit\`a di Firenze, Viale Morgagni 67/A, I--50134 Firenze, Italy
}

\date{}

\maketitle
\thispagestyle{empty}

\begin{abstract}
The \emph{Minkowski product} of unit quaternion sets is introduced and
analyzed, motivated by the desire to characterize the overall variation of 
compounded spatial rotations that result from individual rotations subject 
to known uncertainties in their rotation axes and angles. For a special
type of unit quaternion set, the spherical caps of the 3--sphere $S^3$
in $\rr^4$, closure under the Minkowski product is achieved. Products of 
sets characterized by fixing either the rotation axis or rotation angle, 
and allowing the other to vary over a given domain, are also analyzed.
Two methods for visualizing unit quaternion sets and their Minkowski
products in $\rr^3$ are also discussed, based on stereographic projection
and the Lie algebra formulation. Finally, some general principles for 
identifying Minkowski product boundary points are discussed in the case
of full--dimension set operands.
\end{abstract}

\bigskip
\centerline{{\bf Keywords}: 
Minkowski products, unit quaternions, spatial rotations}
\centerline{3--sphere, stereographic projection, Lie algebra,
boundary evaluation.}

\bigskip
\centerline{e--mail addresses: farouki@ucdavis.edu, gentili@math.unifi.it,}
\centerline{hpmoon@dongguk.edu, stoppato@math.unifi.it.}

\section{Introduction}

The Minkowski sum $A \oplus B$ of two point sets $A,B \in \mathbb{R}^n$ 
is the set of all points generated \cite{minkowski} by the vector sums of 
points chosen independently from those sets, i.e.,
\be
\label{Msum}
A \oplus B \,:=\,
\{\, \a+\b \,:\, \a \in A \;{\rm and}\; \b \in B\, \} \,.
\ee
The Minkowski sum has applications in computer graphics, geometric design, 
image processing, and related fields \cite{ghosh88,hartquist99,kaul,kaul95,kaul92,mditch88,serra}.
The validity of the definition (\ref{Msum}) in $\mathbb{R}^n$ for all
$n \ge 1$ stems from the straightforward extension of the vector sum $\a+\b$ 
to higher--dimensional Euclidean spaces. However, to define a Minkowski 
\emph{product} set
\be
\label{Mproduct}
A \otimes B \,:=\,
\{\, \a\,\b \,:\, \a \in A\;{\rm and}\;\b \in B\, \} \,,
\ee
it is necessary to specify \emph{products} of points in $\mathbb{R}^n$. In 
the case $n=1$, this is simply the real--number product --- the resulting 
algebra of point sets in $\mathbb{R}^1$ is called \emph{interval arithmetic} 
\cite{moore,moore2} and is used to monitor the propagation of uncertainty 
through computations in which the initial operands (and possibly also the 
arithmetic operations) are not precisely determined.

A natural realization of the Minkowski product (\ref{Mproduct}) in $\mathbb
{R}^2$ may be achieved \cite{farouki01} by interpreting the points $\a$ and 
$\b$ as \emph{complex numbers}, with $\a\,\b$ being the usual complex--number 
product. Algorithms to compute Minkowski products of complex--number sets 
have been formulated \cite{farouki00}, and extended to determine Minkowski 
roots and powers \cite{farouki00b,farouki02} of complex sets; to evaluate 
polynomials specified by complex--set coefficients and arguments \cite
{farouki04}; and to solve simple equations expressed in terms of complex--set 
coefficients and unknowns \cite{farouki05}. The Minkowski algebra of complex 
sets introduces rich geometrical structures and has useful applications to
mathematical morphology, geometrical optics, and the stability analysis of 
linear dynamic systems.

In proceeding to higher dimensions, it is natural to consider next the
case of $\mathbb{R}^4$, in which a (non--commutative) ``product of points'' 
may be specified by invoking the quaternion algebra. In this context, the 
study of the Minkowski sum has no obvious and intuitive motivation, but 
the use of unit quaternions to describe spatial rotations provides a 
compelling case for the investigation of Minkowski products in $\mathbb
{R}^4$. Applications in computer animation, robot path planning, 5--axis 
CNC machining, and related fields frequently involve compounded sequences 
of spatial rotations, that are individually subject to certain indeterminacies. 
The set of all possible outcomes of such compounded sequences of indeterminate 
spatial rotations possesses a natural description as the (ordered) Minkowski 
product of unit quaternion sets.

The remainder of this paper is organized as follows. Section~\ref
{sec:quaternions} provides a brief review of the key properties of 
unit quaternions and their interpretation as rotation operators, while 
Section~\ref{sec:sp} discusses the visualization of the set of all 
unit quaternions (namely, the 3--sphere $S^3$ in $\rr^4$) through its 
stereographic projection to $\rr^3$. The concept of a Minkowski product 
of two unit quaternion sets is then introduced in Section~\ref{sec:Mproduct}.
In Sections~\ref{sec:sphsets} and Section~\ref{sec:sproduct}, we focus on 
a specific type of unit quaternion set, the \emph{spherical caps} on 
$S^3$, and show that they exhibit closure under the Minkowski product.
Section~\ref{sec:angles-axes} then considers products of sets more closely 
related to applications, specified by fixing either the rotation axis or 
the angle, and varying the other over a given domain. As an alternative to 
stereographic projection, Section~\ref{sec:lie} describes a Lie algebra 
approach to visualizing unit quaternion sets and their Minkowski products, 
which has the virtue of generating bounded images in $\rr^3$. Finally, 
Section~\ref{sec:productboundaries} discusses general principles for 
identifying boundary points of Minkowski products, giving some necessary 
conditions and a sufficient condition for a product of two points to lie 
on the boundary. Section~\ref{sec:close} summarizes the key results of 
this study and suggests further lines of investigation.

\section{Quaternions and spatial rotations}
\label{sec:quaternions}

Quaternions are ``four--dimensional numbers'' of the form
\[
{\cal A} \,=\, a+a_x\i+a_y\j+a_z\k
\qquad {\rm and} \qquad
{\cal B} \,=\, b+b_x\i+b_y\j+b_z\k \,.
\]
where the elements $\i$, $\j$, $\k$ of the quaternion algebra $\mathbb{H}$ 
obey the multiplication rules $\i^2=\j^2=\k^2=-1$ and $\i\,\j=\k$, $\j\,\k=\i$,
$\k\,\i=\j$. A quaternion ${\cal A}$ may be regarded as comprising scalar 
(real) and vector (imaginary) parts $a=\mbox{scal}({\cal A})$ and $\a=a_x\i
+a_y\j+a_z\k=\mbox{vect}({\cal A})$, and we write ${\cal A}=a+\a$. Real 
numbers and $3$--vectors are subsumed as ``pure scalar'' and ``pure vector'' 
quaternions. The sum and (non--commutative) product of ${\cal A}=a+\a$ and 
${\cal B}=b+\b$ can be expressed using scalar and cross products of vectors 
as
\[
{\cal A}+{\cal B} \,=\, a+b+\a+\b \,, \quad
{\cal A}\,{\cal B} \,=\, ab-\langle\a,\b\rangle+a\,\b+b\,\a+\a\times\b\,.
\]
Every quaternion ${\cal A}=a+\a$ has a \emph{conjugate} ${\cal A}^*=a-\a$,
and a non--negative \emph{magnitude} $|{\cal A}|$ defined by $|{\cal A}|^2
={\cal A}^*\!{\cal A}={\cal A}{\cal A}^*=a^2+|\a|^2$, and one can verify that
$({\cal A}\,{\cal B})^*={\cal B}^*\!{\cal A}^*$ and $|{\cal A}\,{\cal B}|=
|{\cal A}|\,|{\cal B}|$. If $|{\cal A}|\not=0$, the quaternion ${\cal A}$ 
has an \emph{inverse} ${\cal A}^{-1}={\cal A}^*/|{\cal A}|^2$ satisfying
${\cal A}^{-1}\!{\cal A}={\cal A}\,{\cal A}^{-1}=1$, and ${\cal A}^{-1}{\cal B}$ 
and ${\cal B}{\cal A}^{-1}$ specify the \emph{left} and \emph{right} division 
of ${\cal B}$ by ${\cal A}$. We also define an \emph{inner product} $\langle
{\cal A},{\cal B}\,\rangle$ of ${\cal A}$ and ${\cal B}$ (regarded as vectors 
in $\mathbb{R}^4$) by
\be
\label{inner}
\langle{\cal A},{\cal B}\,\rangle \,:=\, a\,b+\langle\a,\b\rangle
\,=\, \mbox{scal}({\cal A}\,{\cal B}^*) \,.
\ee

A \emph{unit} quaternion $\U=u+u_x\i+u_y\j+u_z\k$ satisfies $|\U|=1$, and it 
may be identified with a point on the unit 3--sphere $S^3$ in $\mathbb{R}^4$ 
defined by the equation $u^2+u_x^2+u_y^2+u_z^2=1$. Since a product $\U_1\,\U_2$ 
of two unit quaternions is also a unit quaternion, the points of $S^3$ have 
the structure of a (non--commutative) \emph{group} with respect to the 
quaternion product. Note that the inner product (\ref{inner}) is invariant 
under multiplication of the operands by a unit quaternion,
\be
\label{Uinner}
\langle{\cal A},{\cal B}\,\rangle \,=\, 
\langle\,\U{\cal A},\U{\cal B}\,\rangle \,=\, 
\langle\,{\cal A}\,\U,{\cal B}\,\U\,\rangle \,,
\ee
since such multiplications correspond \cite{duval} to rotations in $\rr^4$.

Any unit quaternion $\U$ may be expressed in the form 
\be
\label{U}
\U \,=\, \cos\half\theta+\sin\half\theta\,\n
\ee
for some angle $\theta\in[\,-\pi,\pi\,]$ and unit vector $\n$. This defines 
a \emph{rotation operator} in $\mathbb{R}^3$ --- for any pure vector $\v$, 
the product $\U\,\v\,\U^*$ also defines a pure vector, corresponding to a 
rotation of $\v$ through angle $\theta$ about an axis defined by $\n$. Note 
that $\U$ and $-\,\U$ define equivalent rotations.

Successive spatial rotations can be replaced by a ``compounded'' rotation
--- the result of consecutively applying rotations $\U_2=\cos\half\theta_2+
\sin\half\theta_2\,\n_2$ and $\U_1=\cos\half\theta_1+\sin\half\theta_1\,\n_1$ 
to $\v$ is $\U_1\,(\,\U_2\v\,\U_2^*\,)\;\U_1^*$, which can be expressed as
$\U\,\v\,\U^*$ with $\U=\U_1\,\U_2$. The non--commutative product captures 
the fact that the result of a sequence of rotations depends upon the \emph
{order} of their application. The rotation angle $\theta$ and axis $\n$ 
of $\U$ are given by
\be
\label{theta}
\cos\half\theta \,=\, \cos\half\theta_1\cos\half\theta_2
 -\sin\half\theta_1\sin\half\theta_2\;\langle\n_1,\n_2\rangle\,, \quad
\ee
\be
\label{n}
\n \,=\, \frac{\sin\half\theta_1\cos\half\theta_2\,\n_1 
\,+\, \cos\half\theta_1\sin\half\theta_2\,\n_2 \,+\,
\sin\half\theta_1\sin\half\theta_2\,\n_1\times\n_2}
{\sin\half\theta} \,.
\ee

\section{Stereographic projection to $\mathbb{R}^3$}
\label{sec:sp}

The set of all unit quaternions occupies the 3--sphere $S^3$ in $\mathbb
{R}^4$. To visualize $S^3$ a \emph{stereographic projection} can be used to 
map it into $\mathbb{R}^3$, just as points on the 2--sphere can be imaged 
onto $\mathbb{R}^2$ to generate a map of the earth's surface. We recall the 
following definition.

\begin{dfn}
Consider the conformal map $\Psi:\hh\setminus\{1\} \to \hh \setminus\{-1\}$ 
defined by
\[
\Psi({\cal Q}):=(1-{\cal Q})^{-1}(1+{\cal Q})=(1+{\cal Q})(1-{\cal Q})^{-1} \,,
\]
and its inverse 
\[
\Phi({\cal Q}):=({\cal Q}+1)^{-1}({\cal Q}-1)=({\cal Q}-1)({\cal Q}+1)^{-1} \,.
\]
The maps $\Psi$ and $\Phi$, or rather their continuous extensions to the 
Alexandroff compactification $\widehat{\hh}:=\hh\cup\{\infty\}$, are called 
\emph{quaternionic Cayley transformations}.
\end{dfn}

As conformal maps, $\Psi$ and $\Phi$ map $n$--spheres and $n$--spaces to 
$n$--spheres and $n$--spaces for $0\leq n \leq 3$. In particular, the next 
result concerns the unit $3$-sphere in $\hh$, denoted $S^3$, and the 
$3$--space of purely imaginary quaternions, denoted $\rr^3$.

\begin{prn}
The stereographic projection from the point $-1$ of $S^3$ to $\rr^3$,
defined by
\be
\label{sp}
u+\u \mapsto \frac{\u}{1+u} \,,
\ee
is the restriction of $\Phi$ to $S^3$.
\end{prn}

\begin{prf}
By direct computation, for all ${\cal Q}=q+\q\in\hh\setminus\{-1\}$ we
have
\ba
\Phi({\cal Q}) \!\! &:=& \!\!
(|{\cal Q}|^2+2q+1)^{-1}({\cal Q}^*+1)({\cal Q}-1)
\nonumber \\
\!\! &=& \!\! (|{\cal Q}|^2+2q+1)^{-1}(|{\cal Q}|^2+2\,\q-1) \,.
\ea
Hence, writing ${\cal Q}=\U=u+\u$ when $|{\cal Q}|=1$, we obtain
\[
\Phi(\U)= \frac{2\,\u}{2+2u} = \frac{\u}{1+u} \,. \quad \QED
\]
\end{prf}

The set of unit quaternions $\U=u+u_x\i+u_y\j+u_z\k \in S^3$ can be 
parameterized in terms of \emph{hyperspherical coordinates} $(\alpha,\beta,
\gamma)$ through the expression
\be
\label{hscoords}
(u,u_x,u_y,u_z) \,=\, (\cos\alpha,\sin\alpha\cos\beta,
\sin\alpha\sin\beta\cos\gamma,\sin\alpha\sin\beta\sin\gamma) \,,
\ee
where $\alpha,\beta\in[\,0,\pi\,]$ and $\gamma\in[\,0,2\pi\,]$. In terms of 
the scalar--vector form (\ref{U}) with $\n=n_x\i+n_y\j+n_z\k$ we have
\[
\theta \,=\, 2\alpha \,, \quad
n_x \,=\, \cos\beta \,, \quad
n_y \,=\, \sin\beta\cos\gamma \,, \quad
n_z \,=\, \sin\beta\sin\gamma \,,
\]
and conversely
\[
\alpha \,=\, \half\theta \,, \quad
\beta \,=\, \arccos n_x \,, \quad
\gamma \,=\, \arctan(n_y,n_z) \,,
\]
where $\arctan(a,b)$ is the angle with cosine $a/\sqrt{a^2+b^2}$ and sine 
$b/\sqrt{a^2+b^2}$.

For each unit quaternion $\U \in S^3$, the point $\v=v_x\i+v_y\j+v_z\k=
\Phi(\U)\in\mathbb{R}^3$ defined by the stereographic projection (\ref{sp}) 
may be identified as the intersection with $\mathbb{R}^3$ of the line in 
$\mathbb{R}^4$ that passes through $-1$ and $\U$. In terms of the 
hyperspherical coordinates (\ref{hscoords}) on $S^3$, we have
\be
\label{sp2}
(v_x,v_y,v_z) \,=\, \tan\half\alpha\,
(\cos\beta,\sin\beta\cos\gamma,\sin\beta\sin\gamma) \,.
\ee
Thus $\v$ may be interpreted as the point with ``ordinary'' spherical 
coordinates $(\beta,\gamma)$ on the 2--sphere in $\mathbb{R}^3$ with radius 
$r=\tan\half\alpha$. The points $\U=1$ and $\U=-1$, corresponding to 
$\alpha=0$ and $\alpha=\pi$, are mapped to the origin of $\mathbb{R}^3$ 
and to infinity, respectively. In terms of the scalar--vector form (\ref{U}) 
of $\U$, the stereographic projection to $\mathbb{R}^3$ becomes
\be
\label{sp3}
\v \,=\, \tan\quarter\theta\;\n \,,
\ee
i.e., $\U$ is mapped to the point identified by the unit vector $\n$ on 
the 2--sphere of radius $r=\tan\quarter\theta$ in $\mathbb{R}^3$. Note that, 
although the unit quaternion $-\,\U=-\cos\half\theta-\sin\half\theta\,\n$ 
specifies a rotation by angle $-\theta$ about $-\,\n$, equivalent to that 
specified by $\U$, it is mapped to the \emph{distinct} point
\[
\v' \,=\, -\cot\quarter\theta\;\n \,.
\]

\section{Quaternionic Minkowski product}
\label{sec:Mproduct}

\begin{dfn}
The \emph{Minkowski product} of two subsets $U,V$ of $\hh$ is 
defined by
\[
U \otimes V \,=\, \{\; \U\,\V \,:\,\U \in U \,,\, \V \in V \;\} \,.
\]
\end{dfn}

If $U,V\subseteq S^3$, then the elements of $U \otimes V$ describe all 
possible compounded rotations generated by a rotation $\V \in V$ followed 
by a rotation $\U \in U$. 

\begin{rmk}
\label{rmk:associative}
As a consequence of the properties of quaternionic multiplication, $\otimes$ 
is an associative but noncommutative operation on the power set of $\hh$.
\end{rmk}

The following Lemma can be verified by direct computation.

\begin{lma}
\label{lemma:differential}
Let 
\begin{align*}
T: \hh \times \hh &\to \hh \\
(\P,\Q)&\mapsto \P\Q\,.
\end{align*}
Then for all $\V\in \hh$ the directional derivatives of $T$ in the directions
$(\V,0)$ and $(0,\V)$ are
\[
\frac{\partial T}{\partial(\V,0)}(\P,\Q)=\V\Q
\quad \mbox{and} \quad
\frac{\partial T}{\partial(0,\V)}(\P,\Q)=\P\V \,.
\]
\end{lma}

\begin{rmk}
Consider two unit quaternions 
\[
\U_1=\cos\half\theta_1+\sin\half\theta_1\,\n_1 , \quad
\U_2=\cos\half\theta_2+\sin\half\theta_2\,\n_2
\]
and the unit quaternions $\U_1+\delta\,\U_1, \U_2+\delta\,\U_2$ that result 
from perturbations $\delta\theta_1,\delta\n_1, \delta\theta_2,\delta\n_2$ to 
the rotation angles and axes of $\U_1,\U_2$. By Lemma~\ref{lemma:differential},
the following equality holds to first order:
\[
(\U_1+\delta\,\U_1)(\U_2+\delta\,\U_2)=
\U_1\,\U_2+\delta\U_1\,\U_2+\U_1\,\delta\U_2 \,.
\]
Thus, to first order in $\delta\theta_1,\delta\n_1$ and $\delta\theta_2,
\delta\n_2$ the product $(\U_1+\delta\,\U_1)(\U_2+\delta\,\U_2)$ is always 
distinct from $\U_1\,\U_2$ --- i.e., the unit quaternion product map 
$S^3 \times S^3 \to S^3$ has no stationary points.
\end{rmk}

We address the problem of determining Minkowski products of different 
types of subsets of $S^3$. For subsets $U,V$ of $S^3$ of full dimension, 
we would ideally like to determine either (i) a ``faithful'' (one--to-one) 
parameterization, over a suitable domain in three parameters, of the 
product set $U \otimes V \subset S^3$; or (ii) a characterization of its 
boundary $\partial(U \otimes V)$ in $S^3$. Problem (ii) is, in general, 
more tractable. We show in Section~\ref{sec:productboundaries} that 
$\partial(U \otimes V) \subseteq\partial U \otimes \partial V$, so (ii) 
amounts to identifying corresponding points $\U\in\partial U$ and $\V\in
\partial V$ that generate (potential) points on the Minkowski product 
boundary $\partial(U \otimes V)$. 

\section{Unit quaternion spherical caps}
\label{sec:sphsets}

In the Minkowski algebra of complex sets \cite{farouki01}, emphasis was 
placed on circular disks as set operands, and it seems natural to extend 
this to the context of unit quaternion sets. For this reason, the first 
class of subsets of $S^3$ that we will study is that of \emph{spherical 
caps}, namely the subsets 
\be
\label{sph}
\{\,\U \in S^3 : |\,\U-\U_0|\leq\rho\,\}
\ee
defined by the intersection of $S^3$ with a $4$--ball that has a prescribed 
radius $\rho$ and unit quaternion $\U_0$ as center. The set (\ref{sph}) 
includes all unit quaternions $\U$ whose distance (measured on $S^3$) from 
$\U_0$ does not exceed $\rho$. In the case of deviations $\delta\,\U$ resulting 
from small perturbations $\delta\theta$ and $\delta\n$ to the rotation angle 
$\theta_0$ and axis $\n_0$ of $\U_0$, it includes all unit quaternions 
satisfying
\[
|\,\delta\,\U\,| \,=\,
\sqrt{\quarter(\delta\theta)^2+\sin^2\!\half\theta_0\,|\delta\n|^2}
\,<\, \rho \,.
\]

\begin{rmk}
\label{rmk:sphericalset}
For $\U_0 \in S^3$ the intersection of $S^3$ with the ball of radius 
$\rho$ and center $\U_0$ in $\hh$ is identical to its intersection with 
a half--space orthogonal to $\U_0$, namely
\[
\{\,\U \in S^3 : |\,\U-\U_0|\leq\rho\,\}=
\{{\,\cal U} \in S^3 : \langle\,\U,\U_0\rangle\geq
1-\half\rho^2\,\} \,,
\]
since $|\,\U-\U_0|^2=2\,(1-\langle\,\U,\U_o\rangle)$. We distinguish three 
cases:
\begin{itemize}
\item for $\rho\geq2$, this set coincides with $S^3$;
\item for $0<\rho<2$ it is a proper subset of $S^3$ --- a 
\emph{spherical cap} --- and its boundary in the topology of $S^3$ 
induced by $\rr^4$ is a $2$--sphere;
\item for $\rho=0$, it is the singleton set $\{\,\U_0\}$.
\end{itemize}
\end{rmk}

Based on the preceding remark, we set $\rho=2\sin\half t$ with $t\in
[\,0,\pi\,]$ so that $1-\half\rho^2=\cos t$, and formulate unit quaternion 
spherical caps as follows.

\begin{dfn}
For $\U_0 \in S^3$ and $t\in[\,0,\pi\,]$, we define
\[
U(\U_0,t):=\{\,\U \in S^3 : 
\langle{\,\cal U},\U_0\rangle\geq \cos t\,\} \,.
\]
and denote its boundary in the topology of $S^3$ as $\partial U(\U_0,t)$.
\end{dfn}

Regarding quaternions $\U \in S^3$ as unit vectors in $\mathbb{R}^4$, we 
can also interpret $U(\U_0,t)$ as the intersection of $S^3$ with the cone 
of vectors $\U$ whose inclinations with $\U_0$ do not exceed $t=2\arcsin
\half\rho$. We note that, if $\U_0=1$, then $\langle{\,\cal U},\U_0\rangle$ 
is just the scalar part of $\U$. For the unit quaternion $\U=\cos\half\theta+
\sin\half\theta\,\n$ corresponding to a rotation through angle $\theta$ about 
the axis $\n$, the condition $|\,\U-1| \le \rho=2\sin\half t$ reduces to
\be
\label{ineq}
\cos\half\theta \,\ge\, 1-\half\rho^2 \,=\, \cos t \,,
\ee
i.e., $|\theta| \leq 2\,t=4\arcsin\half\rho$. The boundary $\partial U(1,t)$ 
corresponds to satisfaction of (\ref{ineq}) with equality, i.e., $|\theta|
=2\,t=4\arcsin\half\rho$. As observed in Remark~\ref{rmk:sphericalset}, 
$\U(1,t)=\{1\}$ if $t=0$ ($\rho=0$), and $\U(1,t)=S^3$ if $t=\pi$ ($\rho=2)$.
 In the case $t=\half\pi$ ($\rho=\sqrt{2}$), it is the set of all unit 
quaternions with any rotation axis $\n$ and rotation angles $\theta\in
[\,-\pi,\pi\,]$.

\begin{rmk}
\label{rmk:specialset}
Setting $\exp(s\,\n)=\cos s+\sin s\,\n$, we have
\[
U(1,t)=\{\,\exp(s\,\n) : 0 \leq s \leq t, |\n|=1\,\} \,.
\]
\end{rmk}

\subsection{Visualization by stereographic projection}

We now consider the stereographic projection of the set $U(\U_0,t)$ 
onto $\rr^3$. 

\begin{prn}
Let $\U_0 \in S^3$ and $t \in (0,\pi)$.
\begin{enumerate}
\item 
If $-1\not\in U(\U_0,t)$ then $\Phi(U(\U_0,t))$ is a closed $3$--ball 
in $\mathbb{R}^3$.
\item 
If $-1\in \partial U(\U_0,t)$ then $\Phi(U(\U_0,t)\setminus
\{-1\})$ is a closed half--space in $\rr^3$.
\item
If $-1\in U(\U_0,t)\setminus\partial U(\U_0,t)$ then $\Phi
(U(\U_0,t)\setminus\{-1\})$ is $\rr^3$ minus an open 3--ball.
\end{enumerate}
Finally, for $t=0$ the image of $U(\U_0,t)$ through $\Phi$ is 
$\{\Phi(\U_0)\}$, and for $t=\pi$ it is $\rr^3$.
\end{prn}

\begin{prf}
Since the statements for the cases $t=0$ and $\pi$ are trivial, we focus on
the case $t \in (0,\pi)$. The set $U(\U_0,t)$ is the intersection of $S^3$ 
with a $4$--ball $B$ centered at $\U_0$ in $\hh$. Thus,
\[
\Phi(U(\U_0,t)) = \Phi(S^3) \cap \Phi(B) = \rr^3 \cap \Phi(B) \,.
\]
where $\Phi(B)$ is a closed subset of $\widehat{\hh}$ that includes $\Phi
(\U_0) \in \rr^3\cup\{\infty\}$ as an interior point. The proposition
then follows from the following observations:
\begin{enumerate}
\item If $-1 \not\in U(\U_0,t)$ then $-1 \not\in B$ and $\Phi(B)$ is 
a closed $4$--ball in $\hh$.
\item If $-1\in \partial U(\U_0,t)$ then $-1 \in \partial B$ and 
$\Phi(B\setminus\{-1\})$ is a closed half--space in $\hh$.
\item If $-1 \in U(\U_0,t)\setminus \partial U(\U_0,t)$ then 
$-1$ is an interior point of $B$ and $\Phi(B\setminus\{-1\})$ is $\hh$ 
minus an open $4$--ball. \QED
\end{enumerate}
\end{prf}

\section{Products of unit quaternion spherical caps} 
\label{sec:sproduct}

We now compute the product $U(\U_0,s) \otimes U(\V_0,t)$. To this end, the 
following remark will be useful.

\begin{rmk}
\label{rmk:singleton}
For all $\U_0,\V_0 \in S^3$ and for all $t\in[\,0,\pi\,]$, we have
\[
U(\U_0,t) \otimes \{\V_0\} =
U(\U_0\V_0,t) = \{\U_0\} \otimes U(\V_0,t) ,
\]
since from (\ref{Uinner}) we note that
$\langle\,\U\V_0^*,\U_0\rangle = \langle\,\U,\U_0\V_0\rangle = 
\langle\,\U_0^*\U,\V_0\rangle$.
\end{rmk}

As a first consequence, it is possible to visualize $U(\U_0,t)=U(1,t)\otimes
\{\U_0\}$ and its boundary $\partial U(\U_0,t)=\partial U(1,t)\otimes\{\U_0\}$ 
as copies of $U(1,t)$ and $\partial U(1,t)$, rotated within $S^3$ so as to 
have center $\U_0$ instead of $1$. We recall that $U(1,t)$ and $\partial 
U(1,t)$ have been described in detail in the first part of Section~\ref
{sec:sphsets}.

Before considering the general product $U(\U_0,s) \otimes U(\V_0,t)$, it 
is instructive to examine a special case.

\begin{lma}
\label{lm:center1}
Let $s,t \in [\,0,\pi\,]$. Then
\[
U(1,s) \otimes U(1,t) = 
\begin{cases}
\,U(1,s+t) & \text{if $s+t \in [\,0,\pi\,]$}, \\ 
\,\qquad S^3 & \text{if $s+t \in [\,\pi,2\pi\,]$}.
\end{cases}
\]
\end{lma}

\begin{prf}
By Remark~\ref{rmk:specialset}, $U(1,s) \otimes U(1,t)$ is the set of 
all products of the form $\exp(a\,\m)\,\exp(b\,\n)$ with $|\m|=|\n|=1$, 
$0\leq a\leq s$, $0\leq b \leq t$. Now, the scalar part of $\exp(a\,\m)\,
\exp(b\,\n)$ is equal to
\[
\cos a\cos b-\sin a\sin b\;\langle\m,\n\rangle,
\]
which is greater than or equal to $\cos a\cos b-\sin a\sin b= \cos(a+b)$. 
If $s+t\in[\,0,\pi\,]$ this bound implies that $U(1,s) \otimes U(1,t) 
\subseteq U(1,s+t)$. If $s+t\in[\,\pi,2\pi\,]$, the bound only implies 
the trivial inclusion $U(1,s) \otimes U(1,t) \subseteq S^3$.

On the other hand, let $\U=u+\u \in S^3$. If we can identify real numbers 
$a$, $b$ with $0\leq a \leq s$ and $0\leq b\leq t$ such that $u=\cos(a+b)$, 
then $\U=\exp((a+b)\,\p)$ for a suitably chosen $\p$ with $|\p|=1$, and we 
conclude that $\U=\exp(a\,\p)\exp(b\,\p) \in U(1,s) \otimes U(1,t)$. If 
$s+t\in[\,0,\pi\,]$, then such $a$ and $b$ exist when $u \geq \cos(s+t)$. 
If $s+t\in[\,\pi,2\pi\,]$, then they exist for all $u\geq-1$. This proves 
that $U(1,s) \otimes U(1,t) \supseteq U(1,s+t)$ in the former case, and 
that $U(1,s) \otimes U(1,t) \supseteq S^3$ in the latter case.
\QED\end{prf}

We are now ready to present the general result for the Minkowski products 
of unit quaternion spherical caps.

\begin{thm}
\label{thm:caps}
Let $\U_0,\V_0 \in S^3$ and $s,t \in [\,0,\pi\,]$. Then
\[
U(\U_0,s) \otimes U(\V_0,t) = 
\begin{cases}
\,U(\U_0\V_0,s+t) & \text{if $s+t \in [\,0,\pi\,]$}, \\ 
\;\;\;\qquad S^3 & \text{if $s+t \in [\,\pi,2\pi\,]$}.
\end{cases}
\]
\end{thm}

\begin{prf}
From Remark~\ref{rmk:singleton}, we have $U(\U_0,s)=\{\U_0\}\otimes U(1,s)$ 
and $U(\V_0,t)=U(1,t)\otimes\{\V_0\}$. Taking into account Remark~\ref
{rmk:associative}, we write
\[
U(\U_0,s) \otimes U(\V_0,t) = 
\{\U_0\}\otimes U(1,s) \otimes U(1,t) \otimes\{\V_0\} \,.
\]
We now apply Lemma~\ref{lm:center1}. If $s+t \in [\,0,\pi\,]$ then
\[
U(\U_0,s) \otimes U(\V_0,t) = \{\U_0\}\otimes U(1,s+t)
\otimes\{\V_0\} \,,
\]
whence
\[
U(\U_0,s) \otimes U(\V_0,t) = U(\U_0\V_0,s+t)
\]
by two further applications of Remark~\ref{rmk:singleton}. However, if 
$s+t\in[\,\pi,2\pi\,]$, then 
\[
U(\U_0,s) \otimes U(\V_0,t) = \{\U_0\}\otimes S^3\otimes
\{\V_0\} = S^3 \,. \quad \QED
\]
\end{prf}

\section{Bounded rotation angles and axes}
\label{sec:angles-axes}

Although unit quaternion spherical caps admit a simple and elegant theory 
for their Minkowski products, they correspond to somewhat complicated and 
non--intuitive relations among the feasible rotation axes and angles. We 
now consider different types of sets, of greater relevance to the physical 
actuators used in robot manipulators, 5--axis milling machines, and related 
contexts.

Specifically, we analyze below the Minkowski products of sets defined by 
(1) fixed rotation axes, and rotation angles that vary over prescribed 
subsets of $[\,-\pi,\pi\,]$; and (2) fixed rotation angles, and rotation 
axes that deviate from prescribed directions by no more than a given angle. 
In case (1) the operand sets are curves and their product is a 2--surface 
in $S^3$, while in case (2) the operand sets are 2--surfaces and their 
product is of full dimension in $S^3$.

\subsection{Fixed rotation axis, bounded angle}

Let us consider the following sets.

\begin{dfn}
For each quaternionic imaginary unit $\c$, and for all $\phi\in\rr$ and 
all $\delta\in[\,0,\pi\,]$,
\[
C(\c,\phi,\delta) := \{\,\exp(s\,\c) : |s-\phi|\leq\delta\,\} \,.
\]
\end{dfn}

$C(\c,\phi,\pi)$ is the great circle in $S^3$ that passes through $1$ and 
$\exp(\phi\,\c)$. For all $\delta\in(0,\pi)$ the set $C(\c,\phi,\delta)$ 
is an arc of this great circle. Finally, $C(\c,\phi,0)$ is the singleton 
$\{\exp(\phi\,\c)\}$.

\begin{prn}
For a fixed quaternionic imaginary unit $\c$, angles $\phi_1,\phi_2\in\rr$,
and ranges $\delta_1,\delta_2\in[\,0,\pi\,]$, the product $C(\c,\phi_1,\delta_1)
\otimes C(\c,\phi_2,\delta_2)$ is a circle, a circular arc, or a singleton:
\[\
C(\c,\phi_1,\delta_1) \otimes C(\c,\phi_2,\delta_2) = 
\begin{cases}
\, C(\c,\phi_1+\phi_2,\delta_1+\delta_2) &
\text{if $\delta_1+\delta_2 \in [\,0,\pi\,]$}, \\
\, C(\c,\phi_1+\phi_2,\pi) &
\text{if $\delta_1+\delta_2 \in [\,\pi,2\pi\,]$}.
\end{cases}
\]
\end{prn}

\begin{prf}
The statement is an immediate consequence of the fact that
\[
\exp(s\,\c)\exp(t\,\c) = \exp((s+t)\c)
\]
for all $ s,t \in \rr$.
\QED\end{prf}

We now study the nature of the product $C(\c_1,\phi_1,\delta_1) \otimes 
C(\c_2,\phi_2,\delta_2)$ in greater detail. To this end, the following remark 
will prove useful.

\begin{rmk}
\label{rmk:circle}
For each quaternionic imaginary unit $\c$, and for all $\phi\in\rr$ and all 
$\delta\in[\,0,\pi\,]$, 
\[
C(\c,0,\delta) \otimes \{\exp(\phi\,\c)\} = C(\c,\phi,\delta) 
= \{\exp(\phi\,\c)\} \otimes C(\c,0,\delta)\,.
\]
\end{rmk}

\begin{thm}
For fixed quaternionic imaginary units $\c_1,\c_2$ with $\c_2\neq\pm\,\c_1$,
the product 
\[
C(\c_1,\phi_1,\delta_1) \otimes C(\c_2,\phi_2,\delta_2)\subset S^3
\]
is an immersed $2$--surface in $\hh$, possibly with boundary. The smallest 
unit quaternion spherical cap $U(\exp(\phi_1\c_1)\exp(\phi_2\c_2),\eta)$ 
that includes this product has $\eta:=\arccos(r)$, where
\begin{equation}
\label{eq:minimum}
r:=\min_{|s|\leq\delta_1,|t|\leq\delta_2} (\cos s\cos t-\sin s\sin t\,
\langle\c_1,\c_2\rangle)\,.
\end{equation}
Moreover, if neither $\delta_1$ nor $\delta_2$ is equal to $\pi$, and at 
least one of them is less than $\half\pi$, then $C(\c_1,\phi_1,\delta_1) 
\otimes C(\c_2,\phi_2,\delta_2) \subset S^3$ is an embedded $2$--surface 
in $\hh$, with boundary. Its boundary consists of four circular arcs which 
have pairwise intersections at the four points 
\[
\exp((\phi_1\pm\delta_1)\c_1)\exp((\phi_2\pm\delta_2)\c_2) \,, \quad 
\exp((\phi_1\pm\delta_1)\c_1)\exp((\phi_2\mp\delta_2)\c_2) \,.
\]
\end{thm}

\begin{prf}
By Remarks~\ref{rmk:singleton} and \ref{rmk:circle}, it suffices to 
consider the case $\phi_1=\phi_2=0$. Consider the surjective map
\begin{align*}
P:[-\delta_1,\delta_1]\times[-\delta_2,\delta_2] &\to C(\c_1,0,\delta_1) 
\otimes C(\c_2,0,\delta_2) \\
(s,t)&\mapsto \exp(s\,\c_1)\exp(t\,\c_2) \,.
\end{align*}
Note that $P$ is non--singular, since the $s$--derivative $\exp(s\,\c_1)\,
\c_1\exp(t\,\c_2)$ and $t$--derivative $\exp(s\,\c_1)\exp(t\,\c_2)\,\c_2=
\exp(s\,\c_1)\,\c_2\exp(t\,\c_2)$ cannot be linearly dependent over $\rr$:
if they were, then $\c_1,\c_2$ would also be linearly dependent, contradicting 
the hypothesis $\c_2\neq\pm\c_1$.

Now let us determine for which $\eta$ the inclusion
\[
C(\c_1,0,\delta_1) \otimes C(\c_2,0,\delta_2) \subseteq U(1,\eta)
\]
holds. The scalar part of the product $\exp(s\,\c_1)\exp(t\,\c_2)$ is
\[
\cos s \cos t - \sin s \sin t\,\langle\c_1,\c_2\rangle \,.
\]
This quantity spans the whole interval $[\,r,1\,]$ with
\[
r:=\min_{|s|\leq\delta_1,|t|\leq\delta_2} 
(\cos s\cos t-\sin s\sin t\,\langle\c_1,\c_2\rangle) \,.
\]

Moreover, $P$ is an embedding if and only if $P$ is injective. The 
equality $\exp(p\,\c_1)\exp(r\,\c_2)=\exp(s\,\c_1)\exp(t\,\c_2)$ holds if, 
and only if, $\exp((p-s)\c_1)=\exp((t-r)\c_2)$, i.e., $p-s,t-r\in\{-2\pi,
0,2\pi\}$ or $p-s,t-r\in\{\pm\pi\}$. Thus, $P$ is injective if and only 
if neither $\delta_1$ nor $\delta_2$ equals $\pi$ and at least one of 
them is less than $\half\pi$. When $P$ is an embedding, the boundary of 
$C(\c_1,0,\delta_1) \otimes C(\c_2,0,\delta_2)$ consists of the four 
circular arcs
\[
[-\delta_1,\delta_1] \to C(\c_1,0,\delta_1) \otimes 
C(\c_2,0,\delta_2)\quad s \mapsto \exp(s\,\c_1)\exp(\pm\delta_2\c_2) \,,
\]
\[
[-\delta_2,\delta_2] \to C(\c_1,0,\delta_1) \otimes C(\c_2,0,\delta_2)
\quad t \mapsto \exp(\pm\delta_1\c_1)\exp(t\,\c_2) \,,
\]
as desired.
\QED
\end{prf}

\begin{prn}
If $\delta_1,\delta_2\in[\,0,\half\pi\,]$ then the smallest 
unit quaternion spherical cap $U(\exp(\phi_1\c_1)\exp(\phi_2\c_2),\eta)$ 
that includes $C(\c_1,\phi_1,\delta_1) \otimes C(\c_2,\phi_2,\delta_2)$ has 
\[
\eta=\arccos\left(\cos\delta_1\cos\delta_2-\sin\delta_1\sin\delta_2
|\langle\c_1,\c_2\rangle|\right)<\delta_1+\delta_2\,.
\]
Furthermore, the boundary of the embedded surface $C(\c_1,\phi_1,\delta_1) 
\otimes C(\c_2,\phi_2,\delta_2)$ intersects the boundary of $U(\exp(\phi_1
\c_1)\exp(\phi_2\c_2),\eta)$ either at the two corners $\exp((\phi_1\pm
\delta_1)\c_1)\exp((\phi_2\pm\delta_2)\c_2)$ or at $\exp((\phi_1\pm\delta_1)
\c_1)\exp((\phi_2\mp\delta_2)\c_2)$.
\end{prn}

\begin{prf}
Let $F(s,t):=\cos s\cos t-\sin s\sin t\,\langle\c_1,\c_2\rangle$. By 
direct computation, the only critical point of $F$ in the interior of 
$[-\delta_1,\delta_1]\times[-\delta_2,\delta_2]$ is its maximum point $(0,0)$. 
Moreover, the restrictions to $(-\delta_1,\delta_1)\times\{\pm\delta_2\}, 
\{\pm\delta_1\}\times(-\delta_2,\delta_2)$ are concave. Hence, the minimum 
of $F$ in $[-\delta_1,\delta_1]\times[-\delta_2,\delta_2]$ is attained at 
one of the four corner points $(\pm\delta_1,\pm\delta_2),(\pm\delta_1,
\mp\delta_2)$. Now
\[
F(\pm\delta_1,\pm\delta_2) = \cos\delta_1\cos\delta_2
- \sin\delta_1\sin\delta_2\,\langle\c_1,\c_2\rangle,
\]
\[
F(\pm\delta_1,\mp\delta_2) = \cos\delta_1\cos\delta_2
+ \sin\delta_1\sin\delta_2\,\langle\c_1,\c_2\rangle,
\]
and consequently
\[
\eta=\arccos\left(\cos\delta_1\cos\delta_2-\sin\delta_1\sin\delta_2
|\langle\c_1,\c_2\rangle|\right)<\delta_1+\delta_2 \,,
\]
since $\cos(\delta_1+\delta_2)=\cos\delta_1\cos\delta_2-\sin\delta_1\sin
\delta_2$.
\QED\end{prf}

\subsection{Bounded rotation axis, fixed angle}

We now consider the following sets.

\begin{dfn}
For each quaternionic imaginary unit $\c$ and for all $\phi\in(0,\pi)$ and 
all $\xi\in[\,0,\pi\,]$, we define
\[
S(\c,\phi,\xi):=\{\, \exp(\phi\,\m) =
\cos\phi+\sin\phi\,\m: \langle \m, \c \rangle \geq \cos \xi \,\}\,.
\]
\end{dfn}

$S(\c,\phi,\pi)$ is the $2$--sphere obtained by intersecting $S^3$ with 
the $3$--space of quaternions whose scalar part is equal to $\cos\phi$. 
For all $\xi\in(0,\pi)$, the set $S(\c,\phi,\xi)$ is a spherical cap of 
that $2$--sphere, whose boundary in the topology of the $2$--sphere is 
a circle:
\[
bS(\c,\phi,\xi)=\{\, \exp(\phi\,\m) = \cos\phi+\sin\phi \,\m: 
\langle \m,\c \rangle = \cos \xi \,\} \,.
\]
Here we introduce the symbol $b$ to distinguish this type of boundary 
from the boundary $\partial$ of the same set in the topology of $S^3$. 
Finally, we note that $S(\c,\phi,0)=\{\,\exp(\phi\,\c)\,\}$.

\begin{prn}
Choose quaternionic imaginary units $\c_1,\c_2$ and let $\phi_1,\phi_2\in
(0,\pi)$ and $\xi_1,\xi_2\in(0,\pi\,]$. Then the rank of the real differential 
of the map
\begin{align*}
S(\c_1,\phi_1,\xi_1) \times S(\c_2,\phi_2,\xi_2) &\to S^3 \\
(\exp(\phi_1\m),\exp(\phi_2\n))&\mapsto \exp(\phi_1\m)\exp(\phi_2\n)
\end{align*}
at a point $(\exp(\phi_1\m),\exp(\phi_2\n))$ is less than $3$ if and only 
if $\m=\pm\n$.
\end{prn}

\begin{prf}
We denote the map by $\sigma$ and fix a point $(\exp(\phi_1\m),
\exp(\phi_2\n)) \in S(\c_1,\phi_1,\xi_1) \times S(\c_2,\phi_2,\xi_2)$. 
By Lemma~\ref{lemma:differential}, we have
\[
\frac{\partial\sigma}{\partial (\v,0)}(\exp(\phi_1\m),\exp(\phi_2\n))
=  \v\exp(\phi_2\n) \,,
\]
\[
\frac{\partial\sigma}{\partial (0,\w)}(\exp(\phi_1\m),\exp(\phi_2\n))
= \exp(\phi_1\m)\,\w \,,
\]
for all $(\v,0)$ and $(0,\w)$ in the tangent space to $S(\c_1,\phi_1,\xi_1) 
\times S(\c_2,\phi_2,\xi_2)$ at the point $(\exp(\phi_1\m),\exp(\phi_2\n))$,
i.e., for all 
\[
\v\in\Pi_\m:=\{\v : \v\perp\m\},\quad \w\in\Pi_\n:=\{\w : \w\perp\n\}.
\]
The image of the differential of $\sigma$ is the sum $\Pi_1+\Pi_2$ of 
the $2$--plane $\Pi_1:=\Pi_\m \exp(\phi_2\n)$ through the origin and the 
$2$--plane $\Pi_2:=\exp(\phi_1\m)\,\Pi_\n$ through the origin. Now, 
$\Pi_1+\Pi_2$ has dimension less than $3$ if and only if $\Pi_1=\Pi_2$,
which is equivalent to 
\[
\exp(-\phi_1\m)\,\Pi_\m=\Pi_\n\exp(-\phi_2\n) \,,
\]
and this in turn is equivalent to $\Pi_\m=\Pi_\n$, i.e., $\m=\pm\n$.
\QED\end{prf}

The previous result immediately implies the following description of the 
product $S(\c_1,\phi_1,\xi_1) \otimes S(\c_2,\phi_2,\xi_2)$ for small 
$\xi_1,\xi_2$.

\begin{crl}
Choose quaternionic imaginary units $\c_1,\c_2$ with $\c_1\neq\pm\c_2$ 
and let $\phi_1,\phi_2\in(0,\pi)$. Then if $\xi_1,\xi_2\in(0,\pi)$ are 
sufficiently small, we have:
\begin{enumerate}
\item the map
\begin{align*}
\sigma:S(\c_1,\phi_1,\xi_1) \times S(\c_2,\phi_2,\xi_2) &\to S^3 \\
(\exp(\phi_1\m),\exp(\phi_2\n)) &\mapsto \exp(\phi_1\m)\exp(\phi_2\n)
\end{align*}
is a submersion;
\item its image $S(\c_1,\phi_1,\xi_1) \otimes S(\c_2,\phi_2,\xi_2)$ is 
(the closure in $S^3$ of) an open subset of $S^3$;
\item the boundary $\partial (S(\c_1,\phi_1,\xi_1) \otimes 
S(\c_2,\phi_2,\xi_2))$ 
is included in 
\[
\left(bS(\c_1,\phi_1,\xi_1) \otimes S(\c_2,\phi_2,\xi_2)\right) 
\cup \left(S(\c_1,\phi_1,\xi_1) \otimes bS(\c_2,\phi_2,\xi_2)\right) \,,
\]
where the two members of the union intersect in the set 
$bS(\c_1,\phi_1,\xi_1) \otimes bS(\c_2,\phi_2,\xi_2)$.
\end{enumerate}
\end{crl}

The next result determines which spherical caps $U(\exp(\phi\,\c),t)$ 
contain $S(\c,\phi,\xi)$ and it provides a rough estimate of which 
$U(\exp(\phi_1\c_1)\exp(\phi_2\c_2),t)$ contain $S(\c_1,\phi_1,\xi_1) 
\otimes S(\c_2,\phi_2,\xi_2)$.

\begin{prn}
\label{prn:estimates}
For each quaternionic imaginary unit $\c$ and for all $\phi\in(0,\pi)$ 
and all $\xi\in[\,0,\pi\,]$, the smallest unit quaternion spherical cap 
$U(\exp(\phi\,\c),t)$ that includes $S(\c,\phi,\xi)$ has $t$ equal to
\[
t(\phi,\xi):=\arccos(\cos^2\phi+\sin^2\phi\cos\xi)\leq\xi\,.
\]
As a consequence, for all $\c_1,\c_2$, all $\phi_1,\phi_2\in(0,\pi)$, and 
all $\xi_1,\xi_2\in[\,0,\pi\,]$, if $T:=t(\phi_1,\xi_1)+t(\phi_2,\xi_2)\in
[\,0,\pi\,]$ then the unit quaternion spherical cap
\[
U(\exp(\phi_1\c_1)\exp(\phi_2\c_2),T)
\]
includes $S(\c_1,\phi_1,\xi_1) \otimes S(\c_2,\phi_2,\xi_2)$.
\end{prn}

\begin{prf}
The first statement is a consequence of the fact that, for $\exp(\phi\,\m)
\in S(\c,\phi,\xi)$, the expression
\[
\langle \exp(\phi\,\m),\exp(\phi\,\c)\rangle =
\cos^2\phi+\sin^2\phi\,\langle\m,\c\rangle
\]
attains its minimum when $\langle\m,\c\rangle=\cos\xi$. The second statement 
follows from the inclusions
\[ 
S(\c_i,\phi_i,\xi_i)\subset U(\exp(\phi_i\,\c_i),t(\phi_i,\xi_i)) \,,
\quad i=1,2
\]
and from Theorem~\ref{thm:caps}.
\QED\end{prf}

\begin{rmk}
\label{rmk:inequalitynotstrict}
The inequality $t(\phi,\xi)\leq\xi$ is strict if and only if $\phi\neq
\half\pi$ and $\xi\neq0$.
\end{rmk}

In some instances, the estimate for the product $S(\c_1,\phi_1,\xi_1)\otimes 
S(\c_2,\phi_2,\xi_2)$ presented in Proposition~\ref{prn:estimates} is sharp.

\begin{exm}
{\rm 
For all quaternionic imaginary units $\c$ and all $\xi\in[\,0,\half\pi\,]$, 
\[
S\left(\c,\half\pi,\xi\right) \otimes S\left(\c,\half\pi,\xi\right) 
\subseteq U(-1,2\xi)
\]
by Proposition~\ref{prn:estimates} and Remark~\ref{rmk:inequalitynotstrict}. 
For a fixed unitary $\v$ orthogonal to $\c$, the product of the two elements
\[
\cos\xi\,\c\pm\sin\xi\,\v
\]
of $S\left(\c,\half\pi,\xi\right)$ has scalar product with $-1$ equal to 
$\cos^2\xi-\sin^2\xi=\cos 2\xi$. Hence, the product of these elements does 
not belong to $U(-1,t)$ if $t<2\xi$.
}
\end{exm}

\section{Lie algebra representation}
\label{sec:lie}

As an alternative to stereographic projection, the Lie algebra $\mathfrak
{so}(3)$ associated with the Lie group ${SO}(3)$ provides a more intuitive 
visualization in $\mathbb{R}^3$ of the Minkowski products of unit quaternion 
sets. In this algebra, spatial rotations are represented by \emph{Euler 
vectors} of the form $\theta\,\n$, where $\theta$ is the rotation angle 
and the unit vector $\n$ defines the rotation axis. With $\theta\in[\,-\pi,
\pi\,]$, any set of spatial rotations lies inside the sphere with center 
at the origin and radius $\pi$ in $\mathbb{R}^3$. This approach avoids 
mapping finite points to infinity, which was the case with stereographic 
projection in Section~\ref{sec:sphsets}.

We recall that the elements of the Lie algebra $\mathfrak{so}(3)$ are 
exactly all skew-symmetric $3\times3$ real matrices. Thus,
\[
\a=a_x\i+a_y\j+a_z\k\ \mapsto\ 
\A=\left[\, \begin{array}{ccc}
0 & -a_z & a_y \\
a_z & 0 & -a_x \\
-a_y & a_x & 0
\end{array} \,\right] 
\]
defines an isomorphism between the vector spaces $\rr^3$ and $\mathfrak{so}
(3)$. It maps any vector $\a$ to the matrix $\A$ associated with the linear 
map $\v\mapsto\a\times\v$. If the inverse isomorphism is denoted by
\[
\mathscr{I}:\mathfrak{so}(3)\to\rr^3 \,,
\] 
then by direct computation one can verify that for all $\A,\B\in\mathfrak
{so}(3)$,
\[
\mathscr{I}([\A,\B]) = \mathscr{I}(\A)\times\mathscr{I}(\B) \,,
\]
where $[\A,\B]=\A\B-\B\A$ is the commutator of $\A$ and $\B$.

Each element of $\mathfrak{so}(3)$ has the form $\theta\,\N$ with $\theta\in
\rr$ and
\[
\N \,=\,
\left[\, \begin{array}{ccc}
0 & -n_z & n_y \\
n_z & 0 & -n_x \\
-n_y & n_x & 0
\end{array} \,\right] \,,
\]
where $\n=n_x\i+n_y\j+n_z\k=\mathscr{I}(\N)$ is a unitary element of 
$\rr^3$.

\begin{rmk}
If $\n=\mathscr{I}(\N)$ is unitary, then by direct computation\footnote
{Here the vector $\n\in\rr^3$ is regarded as a column matrix, subject to the 
matrix product.} we have
\[
\N^2 =
\left[\, \begin{array}{ccc}
n_x^2-1 & n_xn_y & n_xn_z \\
n_xn_y & n_y^2-1 & n_yn_z \\
n_xn_z & n_yn_z & n_z^2-1
\end{array} \,\right] = \n\,\n^T-\I \,,
\]
and
\[
\N^3 = \N\,\n\,\n^T-\N = (\n\times\n)\,\n^T - \N = -\,\N \,.
\]
Consequently, the exponential map 
\[
\exp:\mathfrak{so}(3)\to{SO}(3)
\]
acts on $\theta\,\N$ as follows:
\begin{align*}
\exp(\theta\,\N)
&= \I \,+\, \theta\,\N \,+\,\frac{\theta^2}{2!}\,\N^2 
\,+\, \frac{\theta^3}{3!}\,\N^3 \,+\, \cdots\\
&= \I \,+\, 
\left(\theta-\frac{\theta^3}{3!}+\frac{\theta^5}{5!}-\cdots\right)\N 
\,+\, \left(\frac{\theta^2}{2!}-\frac{\theta^4}{4!}+\frac{\theta^6}{6!}
-\cdots\right)\N^2\\
&= \I \,+\, \sin\theta\;\N \,+\,(1-\cos\theta)\,\N^2 \,.
\end{align*}
\end{rmk}

The explicit form of the matrix $\M\in{SO}(3)$ defining a rotation by angle 
$\theta$ about a unit axis vector $\n=n_x\i+n_y\j+n_z\k$ is well known, e.g.,
\cite[p.~75]{altmann}.

\begin{rmk}
If $\M=\exp(\theta\,\N)$, then $\half\,(\M-\M^T) \,=\, \sin\theta\;\N$. 
Consequently, the logarithmic map $\log: {SO}(3)\to\mathfrak{so}(3)$ acts 
on $\M$ as follows:
\be
\label{log}
\log(\M) \,=\, \frac{\arcsin(\|\mathscr{I}(\A)\|)}{\|\mathscr{I}(\A)\|}\;\A,
\quad \A \,=\, \half\,(\M-\M^T) \,,
\ee
where the range of $\arcsin$ is $[\,-\half\pi,\half\pi\,]$.
\end{rmk}

Multiplication on ${SO}(3)$ induces, through the logarithmic map (\ref{log}),
an operation on $\mathfrak{so}(3)$ defined by the Baker--Campbell--Hausdorff 
(BCH) formula \cite{gilmore,selig}. We give a brief introduction to this 
approach, and illustrate its use in visualizing Minkowski products of unit 
quaternion sets by some examples.

\begin{dfn}[Baker--Campbell--Hausdorff formula]
We define $\BCH:\mathfrak{so}(3)\times\mathfrak{so}(3)\to\mathfrak{so}(3)$ 
as the unique function such that
\[
\log(\M_1\M_2) \,=\, \BCH(\log(\M_1),\log(\M_2))
\]
for all $\M_1,\M_2\in{SO}(3)$. The symbol $\BCH$ is also used to denote the 
function $\rr^3\times\rr^3\to\rr^3$ mapping any pair 
$(\mathscr{I}(\A_1),\mathscr{I}(\A_2))$ to $\mathscr{I}(\BCH(\A_1,\A_2))$.
\end{dfn}

Note that expressions (\ref{theta}) and (\ref{n}) suffice to determine the
rotation angle $\theta$ and axis $\n$ in the relation
\[
\theta\,\n \,=\, \BCH(\theta_1\n_1,\theta_2\n_2) \,.
\]

We now illustrate the Lie algebra approach by several computed examples.

\begin{figure}[htbp]
\includegraphics[width=0.32\textwidth]{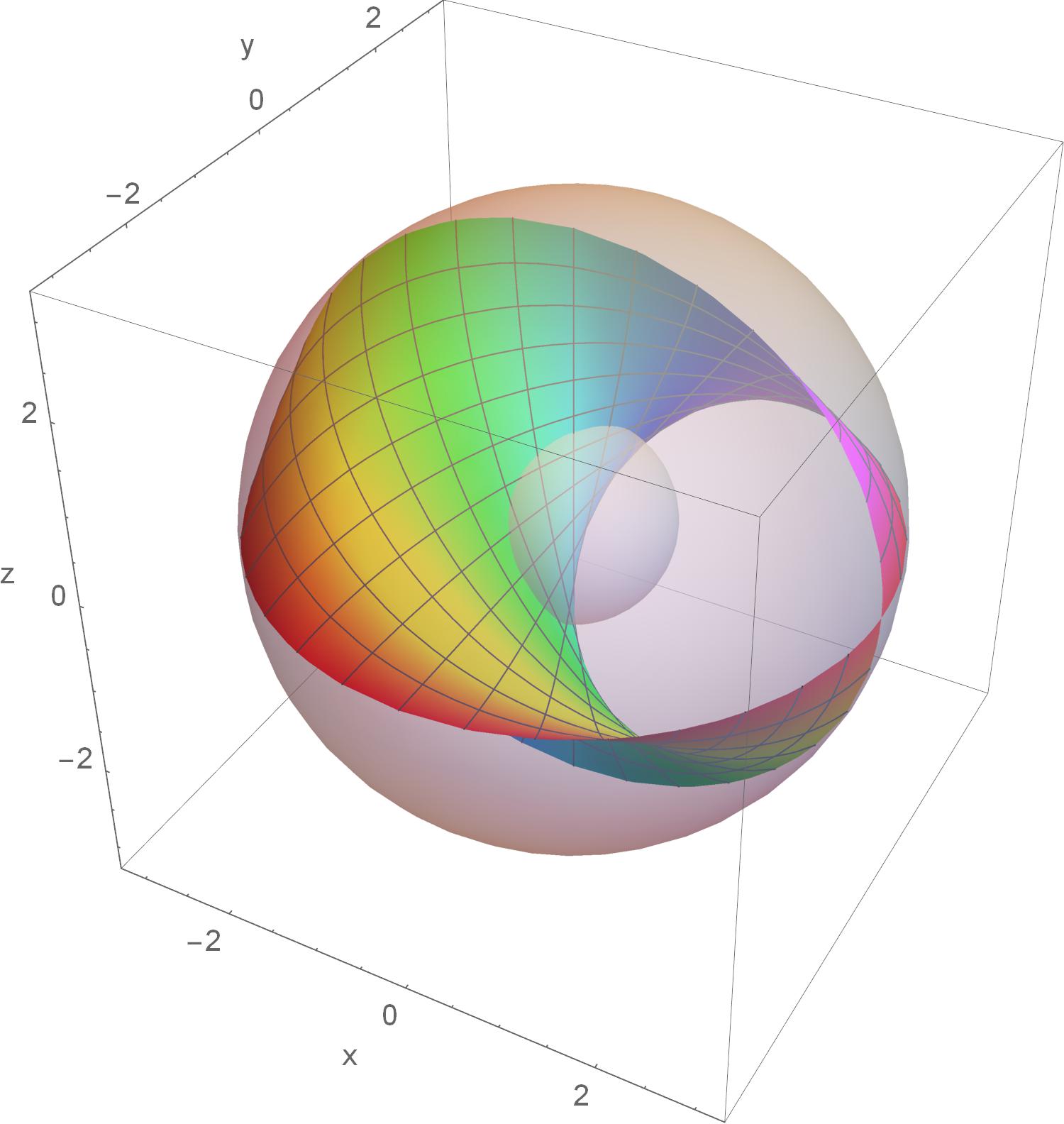}
\includegraphics[width=0.32\textwidth]{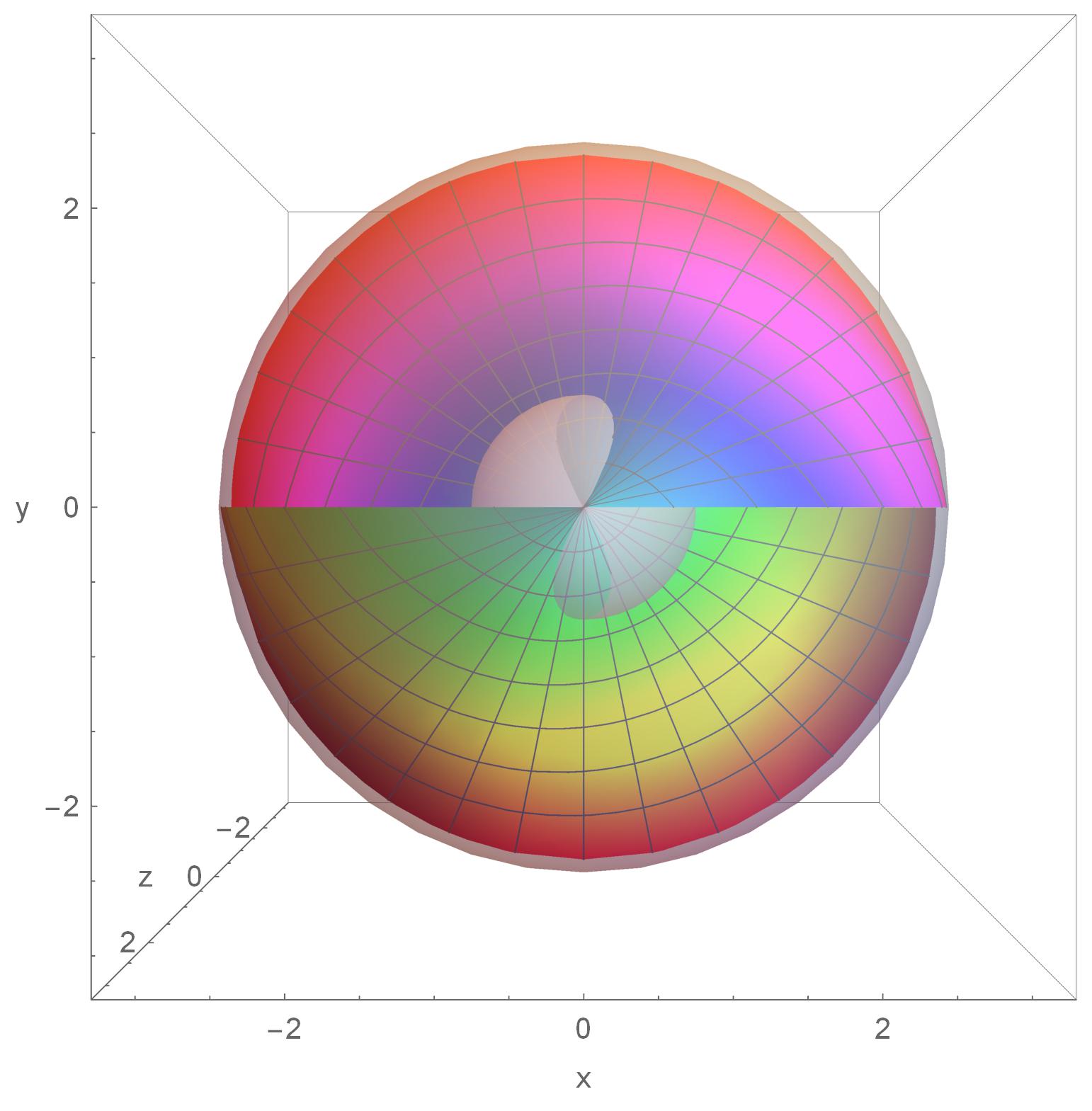}
\includegraphics[width=0.32\textwidth]{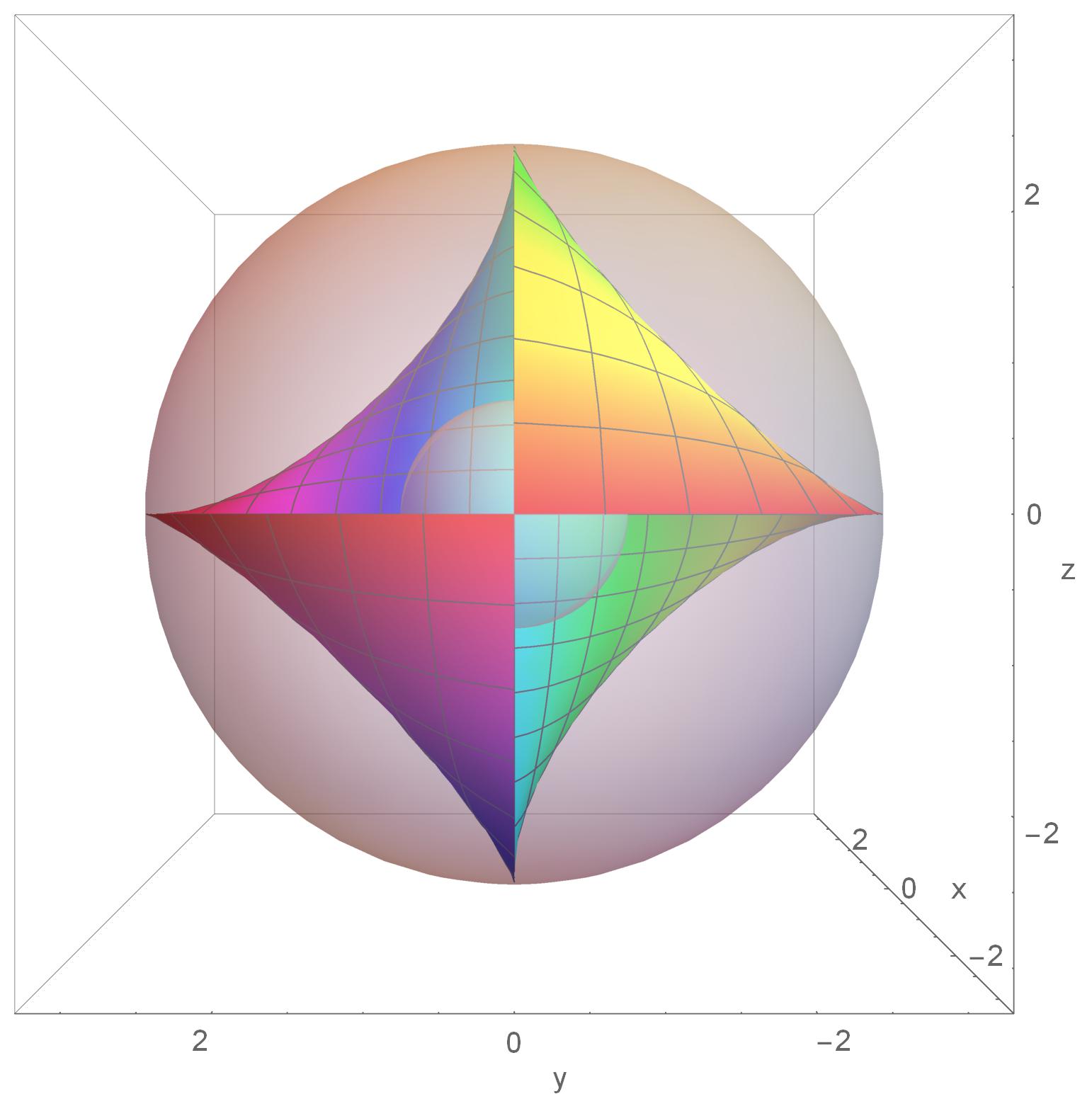}
\caption{A general view (left), top view (center), and left view (right) 
of the Minkowski product of the two great circles $U$ and $V$ specified in 
Example~\ref{exm:example1}.}
\label{fig:example1}
\end{figure}

\begin{exm}
\label{exm:example1}
{\rm We consider (see Section~\ref{sec:angles-axes}) the great circles
\[
U=C(\j,0,\pi) \quad \text{and} \quad V=C(\k,0,\pi)
\]
associated with the subsets ${\bf U}, {\bf V}$ of $SO(3)$ that comprise 
rotations about the $y$ and $z$ axes, respectively. We can visualize the 
Minkowski product $U \otimes V$ as ${\mathscr{I}}(\log{\bf U}{\bf V})=
{\mathrm {BCH}}({\mathscr{I}}(\log{\bf U}),{\mathscr{I}}(\log{\bf V}))
\subset{\mathbb R}^3$. The BCH operands can be parameterized as
\begin{align*}
{\mathscr{I}}(\log{\bf U}) &= \{ t\,\j :\, -2\pi\le t\le2\pi\}, \\
{\mathscr{I}}(\log{\bf V}) &= \{ s\,\k :\, -2\pi\le s\le2\pi\}.
\end{align*}
The set ${\mathscr{I}}(\log{\bf U}{\bf V})$ is a surface lying inside 
the sphere with center 0 and radius $\pi$, as illustrated in Figure~\ref
{fig:example1} (for reference, the spheres of radius 1 and $\pi$ are also
shown). The hue value of the plot corresponds to the parameter $t$.}
\end{exm}

\begin{exm}
\label{exm:example3}
{\rm Consider (see Section~\ref{sec:angles-axes}) the Minkowski product
\[
C(\j,0,\quarter\pi) \otimes S(\k,\eighth\pi,\eighth\pi) \,.
\]
The former set is defined by a fixed rotation axis and variable rotation 
angle, and the latter corresponds to a fixed rotation angle $\theta=
\quarter\pi$ with rotation axes varying in a neighborhood of the $z$ axis. 
The BCH operands are
\begin{align*}
{\mathscr{I}}(\log{\bf U}) &= 
\{ t\,\j :\, -\half\pi\le t\le\half\pi\}, \\
{\mathscr{I}}(\log{\bf V}) &= 
\left\{ \quarter\pi(\cos u\sin v\,\i+\sin u\sin v\,\j+\cos v\,\k) 
:\, -\pi\le u\le \pi,\, 0\le v\le \eighth\pi \right\}.
\end{align*}
The $3$-dimensional set ${\mathscr{I}}(\log{\bf U}{\bf V})$ is the union of 
a one--parameter family of surfaces.}
\end{exm}

\begin{exm}
\label{exm:example5}
{\rm Consider now the Minkowski product
\[
S(\j,\eighth\pi,\eighth\pi)\otimes S(\k,\eighth\pi,\eighth\pi) \,.
\]
These sets have a fixed rotation angle $\theta=\quarter\pi$, and rotation 
axes varying in neighborhoods of the $y$ and $z$ axes, respectively. The 
BCH operands are
\begin{align*}
{\mathscr{I}}(\log{\bf U}) &= 
\left\{ \quarter\pi(\cos s\sin t\,\i+\cos t\,\j+\sin s\sin t\,\k) 
:\, -\pi\le s\le \pi,\, 0\le t\le \eighth\pi \right\}, \\
{\mathscr{I}}(\log{\bf V}) &= 
\left\{ \quarter\pi(\cos u\sin v\,\i+\sin u\sin v\,\j+\cos v\,\k) 
:\, -\pi\le u\le \pi,\, 0\le v\le \eighth\pi \right\}.
\end{align*}
The $3$-dimensional set ${\mathscr{I}}(\log{\bf U}{\bf V})$ can be 
regarded as the union of a two--parameter family of surfaces.}
\end{exm}

\section{Minkowski product boundaries}
\label{sec:productboundaries}

Thus far, we have investigated Minkowski products for specific 1--dimensional, 
2--dimensional, and 3--dimensional subsets of $S^3$. We now consider properties 
of the Minkowski product $U \otimes V$ valid for any subsets $U,V$ of $S^3$ 
that have interior points $\U,\V$ in the topology of $S^3$ induced by $\rr^4$. 
We note that, in this setting, $\U$ and $\V$ are interior points if and only 
if $U$ and $V$ include unit quaternion spherical caps $U(\U,\delta)$ and 
$V(\V,\epsilon)$ with $\delta>0$ and $\epsilon>0$.

\begin{thm}
Let $U,V\subseteq S^3$ and let $\U\in U$, $\V\in V$. Then $\U\V$ is an 
interior point of $U\otimes V$ if $\U$ is an interior point of $U$, or $\V$ 
is an interior point of $V$.
\end{thm}

\prf
When $\U$ is an interior point of $U$, there exists a $\delta>0$ such that 
$U(\U,\delta)\subseteq U$. Then by Remark~\ref{rmk:singleton}, 
\[
U(\U\V,\delta) = U(\U,\delta)\otimes\{\V\} \,,
\]
whence $U(\U\V,\delta)$ is included in  $U\otimes V$, so $\U\V$ is an 
interior point of $U\otimes V$. The case when $\V$ is an interior point of 
$V$ is treated in the same fashion.
\QED

\begin{crl}\label{crl:boundary}
Let $U,V$ be subsets of $S^3$ and let $\partial U,\partial V$ be their 
boundaries in the topology of $S^3$. Then the boundary $\partial(U \otimes V)$ 
of $U\otimes V$  in the topology of $S^3$ is included in $\partial U \otimes 
\partial V$.
\end{crl}

The requirement that $\U \in \partial U$, $\V \in \partial V$ is a 
\emph{necessary} but not \emph{sufficient} condition for $\U\,\V \in 
\partial(U \otimes V)$. In general, products of points on $\partial U$ 
and $\partial V$ may generate interior points of $U \otimes V$, and it 
can be difficult to identify only those pairs of boundary points such 
that $\U\,\V\in\partial(U \otimes V)$:
\begin{itemize}
\item $U \otimes V$ may cover all of $S^3$ (and thus have no boundary) 
even in cases where $U$ and $V$ are proper subsets of $S^3$;
\item it may happen for $\U\in\partial U, 
\V\in\partial V$ that
\[\U\,\V=\U'\,\V'\,,\]
where $\U'$ is an interior point of $U$ or $\V'$ is an interior point of $V$;
\item if $\partial U,\partial V, \partial (U \otimes V)$ are $2$-surfaces in 
$S^3$ then we cannot expect $\partial U \otimes \partial V$ to coincide 
with $\partial (U \otimes V)$ by dimensional considerations.
\end{itemize}

All three phenomena can be observed in the next example:

\begin{exm}
{\rm By Theorem~\ref{thm:caps}, the equality $U(1,s)\otimes U(1,s)=
U(1,2s)$ holds for all $s\in(0,\half\pi]$. Now:
\begin{itemize}
\item if $s=\half\pi$, then $U(1,2s)=S^3$ has no boundary;
\item the product $1$ of the boundary points $\exp(s\i)$ and $\exp(-s\i)$ 
equals the product of the interior points $\exp(\half s\i)$ and 
$\exp(-\half s\i)$;
\item if $s<\half\pi$ then, by direct computation, the product of two 
boundary points $\exp(s\m),\exp(s\n)$ of $U(1,s)$ belongs to $\partial U(1,2s)$ 
if, and only if, $\m=\n$; therefore,
\[(\partial U \otimes \partial U) \setminus \{\U^2 \,:\, \U\in\partial U\}\]
is included in the interior of $U(1,2s)$.
\end{itemize}
}
\end{exm}

For the case of Minkowski sums $S_1 \oplus S_2$ of point sets $S_1,S_2$ 
in $\mathbb{R}^2$ or $\mathbb{R}^3$, a necessary condition for the sum of 
boundary points $\p_1 \in \partial S_1$ and $\p_2 \in \partial S_2$ to 
belong to $\partial(S_1 \oplus S_2)$ is well--known: namely, the normals 
$\n_1,\n_2$ to $\partial S_1,\partial S_2$ at $\p_1,\p_2$ must be linearly 
dependent --- or, equivalently, the tangent spaces to $\partial S_1$ and
$\partial S_2$ at $\p_1$ and $\p_2$ must be identical. This principle was 
extended to Minkowksi products of point sets in $\mathbb{R}^2$ by interpreting 
points as complex numbers $\p_1=x_1+\ri y_1$, $\p_2=x_2+\ri y_2$ and invoking 
the complex logarithm to transform Minkowski products into Minkowski sums 
\cite{farouki00,farouki01}.

The situation with the Minkowski products of unit quaternion sets is 
more subtle, because the map $\log:{SO}(3)\to\mathfrak{so}(3)$ is not a 
homomorphism from the multiplicative group ${SO}(3)$ to the additive 
group $\mathfrak{so}(3)$, i.e., in general
\[
\log(\M_1\M_2)=\BCH(\log(\M_1),\log(\M_2))\neq\log(\M_1)+\log(\M_2) \,.
\]

Nevertheless, we can find a necessary condition using a different strategy. 
We begin by considering the case of spherical caps, before proceeding to 
state and prove a general result.

\begin{lma}
\label{lemma:capproductboundary}
Let $\U_0,\V_0\in S^3$ and let $s_0,t_0$ be such that $0<s_0\leq t_0<\half
\pi$. Suppose that the boundaries in $S^3$ of the spherical caps $U(\U_0,s_0)$ 
and $U(\V_0,t_0)$ intersect at $1$. Then $1$ belongs to the boundary of the 
Minkowski product $U(\U_0,s_0)\otimes U(\V_0,t_0)=U(\U_0\V_0,s_0+t_0)$ in 
$S^3$ if, and only if, there exists an imaginary unit $\c$ such that 
\[
\U_0=\exp(s_0\c),\quad\V_0=\exp(t_0\c) \,.
\]
This happens if, and only if, $U(\U_0,s_0)$ is included in $U(\V_0,t_0)$ and 
the boundaries $\partial U(\U_0,s_0),\partial U(\V_0,t_0)$ are tangent at $1$.
\end{lma}

\prf
By hypothesis, the scalar part of $\U_0$ equals $\cos s_0$ and the scalar 
part of $\V_0$ equals $\cos t_0$. In this case, the scalar part of $\U_0\V_0$ 
equals $\cos(s_0+t_0)=\cos s_0\cos t_0-\sin s_0\sin t_0$ if, and only if, 
there exists an imaginary unit $\c$ such that $\U_0=\exp(s_0\c)$ and $\V_0
=\exp(t_0\c)$. This happens if, and only if, $\U_0$ is included in the 
(smaller) arc of the great circle in $S^3$ with endpoints $1$ and $\V_0$. 
This is the same as asking for $U(\U_0,s_0)$ to be included in $U(\V_0,t_0)$ 
and for $\partial U(\U_0,s_0),\partial U(\V_0,t_0)$ to be tangent at $1$.
\QED

\begin{dfn}
\label{def:tame}
Let $U\subset S^3$. Then the boundary of $U$ in $S^3$ is \emph{tame} if, 
for every point $\U_0 \in \partial U$, the following conditions hold:
\begin{enumerate}
\item[(1)] $U$ includes a spherical cap $U(\U,s)$ with $s>0$, such that 
$U(\U,s)\cap\partial U=\{\U_0\}$;
\item[(2)] at $\U_0$, the set of all tangent vectors to curves in $U$ 
through $\U_0$ spans a $2$--plane, which we denote by $T_{\U_0}\partial U$ 
and call the tangent plane to $\partial U$ at $\U_0$, as usual.
\end{enumerate}
\end{dfn}

In the situation described in the previous definition, $T_{\U_0}\partial U$ 
coincides with the tangent plane to the $2$--sphere $\partial U(\U,s)$ at 
$\U_0$.

\begin{rmk}
If $U$ is the closure in $S^3$ of an open connected subset of $S^3$ and its 
boundary $\partial U$ is a smooth surface, then $\partial U$ is tame.
\end{rmk}

\begin{thm}
\label{thm:UV}
Let $U$, $V$ be proper subsets of $S^3$ and assume that $\partial U$,
$\partial V$ are both tame. Let $\U\in\partial U$ and $\V\in\partial V$. 
If $\U\V$ belongs to $\partial(U \otimes V)$, then 
\[
\U^*\,(T_{\U}(\partial U))= (T_{\V}(\partial V))\,\V^*.
\]
Furthermore, for each spherical cap $U(\U_0,s)\subseteq U$ with $s\in
(0,\half\pi)$ whose boundary is tangent to $\partial U$ at $\U$, and each 
spherical cap $U(\V_0,t)\subseteq V$ with $t\in(0,\half\pi)$ whose boundary 
is tangent to $\partial V$ at $\V$, we have
\[
\U^*\,U(\U_0,s)\subseteq U(\V_0,t)\,\V^*
\quad \mbox{or} \quad
\U^*\,U(\U_0,s)\supseteq U(\V_0,t)\,\V^* \,.
\]
\end{thm}

\prf
Fix any spherical cap $U(\U_0,s)\subseteq U$ whose boundary is tangent 
to $\partial U$ at $\U$, and any spherical cap $U(\V_0,t)\subseteq V$
whose boundary is tangent to $\partial V$ at $\V$. We rotate $U$ to 
$U_1:=\U^*U$, and $V$ to $V_1:=V\V^*$ in $S^3$. After these rotations:
\begin{itemize}
\item the point $1$ belongs to both $\partial U_1$ and $\partial V_1$;
\item $\U^*\,U(\U_0,s)$ is a spherical cap $U(\P,s)\subseteq U_1$ tangent 
to $\partial U_1$ at $1$;
\item $U(\V_0,t)\,\V^*$ is a spherical cap $U(\Q,t)\subseteq V_1$ tangent 
to $\partial V_1$ at $1$;
\item $T_1\partial U_1 = \U^*\,(T_{\U}(\partial U)),\quad T_1\partial V_1 
= (T_{\V}(\partial V))\,\V^*$.
\end{itemize}
Finally, since $U_1 \otimes V_1$ is obtained from $U \otimes V$ by means 
of the rotation $\W\mapsto\U^*\W\V^*$, we conclude that $\U\V\in\partial
(U \otimes V)$ is equivalent to $1\in\partial(U_1 \otimes V_1)$.

We now show that $1\in\partial(U_1 \otimes V_1)$ implies both the equality 
$T_1\partial U_1 =T_1\partial V_1$ and one of the inclusions $U(\P,s)\subseteq 
U(\Q,t), U(\P,s)\supseteq U(\Q,t)$. We claim that
\[
1\in\partial(U(\P,s)\otimes U(\Q,t)) = \partial U(\P\Q,s+t) \,.
\]
Indeed, if $1$ were an interior point of $U(\P\Q,s+t)\subseteq U_1 \otimes 
V_1$ then it would be an interior point of $U_1 \otimes V_1$, contradicting 
our hypothesis. Lemma~\ref{lemma:capproductboundary} allows us to deduce 
that one of $U(\P,s)$ and $U(\Q,t)$ is included in the other, and that 
their boundaries are tangent at $1$. As a consequence,
\[
T_1\partial U_1=T_1\partial U(\P,s)=T_1\partial U(\Q,t)=T_1\partial V_1 \,,
\]
as desired.
\QED

The previous theorem immediately implies the next corollary: a criterion to 
identify parts of $\partial U\otimes\partial V$ that are \emph{not} included 
in $\partial(U\otimes V)$.

\begin{crl}
Let $U,V$ be proper subsets of $S^3$ with tame boundaries, and let $\U\in
\partial U$, $\V\in\partial V$. Suppose two spherical caps 
$U(\U_0,s)$, $U(\U'_0,s')$ exist, such that
\begin{enumerate}
\item neither cap includes the other;
\item both caps are included in $U$;
\item the boundaries of both caps are tangent to $\partial U$ at $\U$.
\end{enumerate}
Then
\[
\U\V\not\in\partial(U\otimes V)
\quad \mbox{and} \quad
\V\U\not\in\partial(V\otimes U) \,.
\]
In other words, $\partial(U\otimes V)$ does not include any point of 
$\{\U\}\otimes V$, and $\partial(V\otimes U)$ does not include any point 
of $V\otimes\{\U\}$.
\end{crl}

\begin{exm}
{\rm Let $s\in(0,\half\pi)$ and set 
\[
U:=U(\exp(s\,\i),s)\cup U(\exp(-s\,\i),s) \,.
\]
Consider the boundary point $1$ of $U$. Then as a consequence of 
the previous corollary, for each proper subset $V$ of $S^3$ with tame boundary, 
the boundary $\partial(U\otimes V)$ does not include any point of $V$.}
\end{exm}

Theorem~\ref{thm:UV} inspires our last result, which identifies a sufficient 
condition for a product $\U\V$ of points $\U\in\partial U$ and $\V\in\partial 
V$ to belong to $\partial(U\otimes V)$.

\begin{thm}
Let $\U\in\partial U$, $\V\in\partial V$ for proper subsets $U$, $V$ of 
$S^3$. Assume that there exist $s\in(0,\half\pi)$ and a spherical cap 
$U(\P,s)$ with $1\in\partial U(\P,s)$ and with
\[\
\U^*U\subseteq U(\P,s)\supseteq V\V^* \,.
\]
Then
\[\U\V\in\partial(U\otimes V)\,.\]
\end{thm}

\prf
Suppose, by contradiction, that $\U\V$ is an interior point of 
\[
U\otimes V\subseteq (\U\,U(\P,s))\otimes(U(\P,s)\,\V)=U(\U\P^2\V,2s) \,.
\]
Then $1$ would be an interior point of 
\[
U(\P^2,2s)=U(\P,s)\otimes U(\P,s) \,,
\]
and since $1\in\partial U(\P,s)$, this would contradict Lemma~\ref
{lemma:capproductboundary}.
\QED

We mention that, in order to apply this sufficient condition, it is not 
necessary to assume that the tangent planes $T_{\U}(\partial U)$, $T_{\V}
(\partial V)$ are well--defined. However, if they are well-defined, then
\[
\U^*T_{\U}(\partial U)=T_1 U(\P,s)=T_{\V}(\partial V)\V^* \,.
\]
We also point out that, in the statement, the assumption $s<\half\pi$ is 
essential. This fact is illustrated in the next examples.

\begin{exm}
{\rm If $\P=\i$ and $s=\half\pi$ then, although $1$ is a boundary point 
of $U(\P,s)$, the Minkowski product
\[
U(\P,s)\otimes U(\P,s)=U(-1,\pi)=S^3
\]
admits no boundary points.}
\end{exm}

\begin{exm}
{\rm Fix $s\in(\half\pi,\pi)$ and let $\P:=\exp(s\i)$, so that 
$1\in\partial U(\P,s)$. Set
\[
U:=U(\P,s)\cap U(1,\eighth\pi)
\]
and consider the boundary point $1$ of $U$. Although $U\subseteq U(\P,s)$, 
the point $1$ is an interior point of $U\otimes U$. Indeed, it can be obtained 
not only as the product of $1$ with itself, but also as the product of the 
two points $\exp(\textstyle\frac1{16}\displaystyle\pi\j)$ and 
$\exp(-\textstyle\frac1{16}\displaystyle\pi\j)$. These two points are interior 
points of $U$ because they are interior points of $U(1,\eighth\pi)$ and because 
\[\langle\exp(\pm\textstyle\frac1{16}\displaystyle\pi\j),\exp(s\i)\rangle=
\cos(\textstyle\frac1{16}\displaystyle\pi)\cos(s)>\cos(s)\,.\]
}
\end{exm}

\section{Closure}
\label{sec:close}

The characterization of the uncertainties in a compounded spatial rotation, 
arising from the ordered product of a sequence of individual rotations 
subject to prescribed uncertainties in their rotation angles and axes, is 
of fundamental interest in diverse contexts. Unit quaternion sets offer 
compact and intuitive representations for families of spatial rotations, 
and their Minkowski products describe the outcomes of all combinations of 
the individual rotations.

Whereas Minkowski sums have been extensively studied, and methods for 
computing them readily generalize to higher dimensions, Minkowski products 
have received less attention. Using complex--number multiplication to define 
products of points in $\rr^2$, algorithms for the Minkowski products of 
planar sets are based upon invoking the complex logarithm to transform 
Minkowski products into Minkowski sums. However, the situation with unit 
quaternion sets is much more challenging, since (i) they reside in a 
non--Euclidean space, the 3--sphere $S^3$; and (ii) their product is 
non--commutative.

This preliminary study of the Minkowski products of unit quaternion sets
describes some basic results and guiding principles for their systematic 
study. A family of unit quaternion sets that is closed under the Minkowski
product --- the spherical caps of $S^3$ --- was identified, and some key
results concerning sets defined by fixing either the rotation axis or 
rotation angle, and allowing the other to vary over a given domain, were 
also developed. To help visualize unit quaternion sets and their Minkowski 
products, mappings to $\rr^3$ based upon stereographic projection and the 
associated Lie algebra were proposed. Finally, general principles to identify 
boundary points of Minkowski products, for full--dimension operand sets with 
smooth boundaries, were analyzed.

The results presented herein serve to introduce the problem of computing
Minkowksi products of unit quaternions, to establish some basic foundations, 
and to identify some possibilities and difficulties it entails. Through its 
depth and practical importance, this problem offers scope for much further 
study.

\raggedright

\subsection*{Acknowledgements}

{\small
The second and fourth authors are partly supported by the Italian Ministry 
of Education (MIUR) through Finanziamento Premiale FOE 2014 ``Splines for 
accUrate NumeRics: adaptIve models for Simulation Environments'' and by 
Istituto Nazionale di Alta Matematica (INdAM) through Gruppo Nazionale per 
le Strutture Algebriche, Geometriche e le loro Applicazioni (GNSAGA).
}


\begin{thebibliography}{99}

\bibitem{altmann}
S.~L.~Altmann (1986), 
\emph{Rotations, Quaternions, and Double Groups}, 
Dover Publications (reprint), Mineola, NY.

\bibitem{duval}
P.~Du~Val (1964), 
\emph{Homographies, Quaternions, and Rotations}, 
Clarendon Press, Oxford.

\bibitem{farouki00b}
R.~T.~Farouki, W.~Gu, and H.~P.~Moon (2000),
Minkowski roots of complex sets,
\emph{Geometric Modeling and Processing 2000}, IEEE Computer Society Press, 
287--300.

\bibitem{farouki04}
R.~T.~Farouki and C.~Y.~Han (2004), 
Computation of Minkowski values of polynomials over complex sets,
\emph{Numerical Algorithms} {\bf 36}, 13--29.

\bibitem{farouki05}
R.~T.~Farouki and C.~Y.~Han (2005), 
Solution of elementary equations in the Minkowski geometric algebra
of complex sets,
\emph{Advances in Computational Mathematics} {\bf 22}, 301--323.

\bibitem{farouki00}
R.~T.~Farouki, H.~P.~Moon, and B.~Ravani (2000), 
Algorithms for Minkowski products and implicitly--defined complex sets,
\emph{Advances in Computational Mathematics} {\bf 13}, 199--229.
 
\bibitem{farouki01}
R.~T.~Farouki, H.~P.~Moon, and B.~Ravani (2001), Minkowski geometric
algebra of complex sets, \emph{Geometriae Dedicata} {\bf 85}, 283--315.
 
\bibitem{farouki02}
R.~T.~Farouki and H.~Pottmann (2002), 
Exact Minkowski products of $N$ complex disks,
\emph{Reliable Computing} {\bf 8}, 43--66.

\bibitem{ghosh88}
P.~K.~Ghosh (1988), A mathematical model for shape description using
Minkowski operators, \emph{Computer Vision, Graphics, and Image Processing} 
{\bf 44}, 239--269.

\bibitem{gilmore}
R.~Gilmore (2008), \emph{Lie Groups, Physics, and Geometry: An Introduction
for Physicists, Engineers, and Chemists},
Cambridge University Press, Cambridge.

\bibitem{hartquist99} 
E.~E.~Hartquist, J.~P.~Menon, K.~Suresh, H.~B.~Voelcker, and J.~Zagajac
(1999), A computing strategy for applications involving offsets, sweeps, 
and Minkowski operators, \emph{Comput.\ Aided Design} {\bf 31}, 175--183.

\bibitem{kaul}
A.~Kaul (1993), Computing Minkowski sums, PhD Thesis, Columbia University.
 
\bibitem{kaul95}
A.~Kaul and R.~T.~Farouki (1995), Computing Minkowski sums of plane curves,
\emph{International Journal of Computational Geometry and Applications} 
{\bf 5}, 413--432.
 
\bibitem{kaul92}
A.~Kaul and J.~R.~Rossignac (1992), Solid interpolating deformations:
Construction and animation of PIP, \emph{Computers and Graphics} {\bf 16}, 
107--115.

\bibitem{mditch88} 
A.~E.~Middleditch (1988), Applications of a vector sum operator, 
\emph{Comput.\ Aided Design} {\bf 20}, 183--188.
 
\bibitem{minkowski}
H.~Minkowski (1903), Volumen und Oberfl\"ache, \emph{Mathematische Annalen} 
{\bf 57}, 447--495.

\bibitem{moore} 
R.~E.~Moore (1966), \emph{Interval Analysis}, Prentice--Hall, Englewood 
Cliffs, NJ.

\bibitem{moore2} 
R.~E.~Moore (1979), \emph{Methods and Applications of Interval Analysis}, 
SIAM, Philadelphia.

\bibitem{selig}
J.~M.~Selig (1996), \emph{Geometrical Methods in Robotics}, Springer, 
New York.
 
\bibitem{serra}
J.~Serra (1982), \emph{Image Analysis and Mathematical Morphology}, Academic 
Press, London.
 
\end{thebibliography}
\end{document}